\documentclass[12pt]{amsart} 
\usepackage{amsmath}
\usepackage{amsfonts}
\usepackage{amsthm}
\usepackage{graphicx}
\usepackage{supertabular}
\usepackage{longtable}
\usepackage{ifthen}
\usepackage{ams}
\usepackage{amssymb}
\DeclareMathOperator{\im}{im} \subjclass[2000]{Primary: 20J06;
Secondary: 46L05, 43A15}

\author{Anthony J. Narkawicz (Virginia Tech)}

\address{502 I Jackson Street\\
Blacksburg, VA 24060, U.S.A}\email{anarkawi@vt.edu}
\title{The First Cohomology Group $H^1(G,M)$}
\begin{document}
\maketitle

\newtheorem{lemma}{Lemma}
\newtheorem{theorem}{Theorem}
\newtheorem{corollary}{Corollary}
\newtheorem{definition}{Definition}
\newtheorem{construction}{Construction}
\newtheorem{conjecture}{Conjecture}

\noindent{\bf Abstract.} This paper characterizes the first
cohomology group $H^1(G,M)$ where $M$ is a Banach space (with norm
$||$ $||_M$) that is also a left $\C G$ module such that the
elements of $G$ act on $M$ as continuous $\C$-linear
transformations. We study this group for $G$ an infinite, finitely
generated group. Of particular interest are the implications of
the vanishing of the group $H^1(G,M)$. The first result is that
$H^1(G,\C G)$ imbeds in $H^1(G,M)$ whenever $\C G \subset M
\subset L^p (G)$ for some $p \in \N.$ This is an unpublished
result and shows immediately that if $H^1(G,M) = 0,$ then $G$ can
have only 1 end. Secondly (also a new result), we show that
$H^1(G,M)$ is not Hausdorff if and only if there exist $f_i \in M$
with norm 1 ($||f_i||_M = 1$) for all $i$ with the property that
$||gf_i - f_i||_M \longrightarrow 0$ as $i \longrightarrow \infty$
for every $g \in G.$ This is then used to show that if $M$ and
$||$ $||_M$ satisfy certain properties and if $G$ satisfies a
``strong F$\o$lner condition," then $H^1(G,M)$ is not Hausdorff.
For the second half of this paper, we give several applications of
these last two theorems focusing on the group $G = \Z^n.$

\vspace{.1 in}

\noindent {\bf 1. Introduction}

\vspace{.1 in}

\noindent Motivation for this paper comes from two papers, one by
Mohammed E.B. Bekka and Alain Valette ([3]) and the other, an
expository paper, by Edward G. Effros ([2]). The first paper
examines the group $H^1(G,L^2(G))$ and focuses on the implications
of the vanishing of this group. It shows the following:

\vspace{.05 in}

(1) $H^1(G,L^2(G))$ is Hausdorff if and only if $G$ is
non-amenable.

(2) The $G$-module imbedding $\C G \to L^p(G)$ induces an
imbedding of

\noindent$H^1(G,\C G)$ into $H^1(G,L^p(G)), p\geq 1.$

(3) If $H^1(G,L^2(G)) = 0,$ then $G$ is non-amenable with just one
end.

\vspace{.05 in}

The first of these results is due in part to the following result
by A. Guichardet ([1],Corollary 2.3 of Chapter III):
$H^1(G,L^p(G))$ is not Hausdorff if and only if there exists a
sequence $e_n$ in $L^p(G)$ such that $||e_n||_p = 1$ for all $n$
with the property that $||ge_n - e_n||_p \longrightarrow 0$ for
all $g \in G$. We show that $L^p(G)$ (with norm $||$ $||_p$) may
be replaced by any Banach space $M$ (with norm $||$ $||_M$) that
is a $\C G$ module and has the property that the elements of $G$
act on $M$ as continuous $\C$-linear transformations. This is not
trivial from Guichardet's theorem. In fact, the topology on
$H^1(G,L^p(G))$ (induced by the $L^p$-norm topology on $L^p(G)$)
is entirely different from the topology on $H^1(G,M)$ which is
induced by the norm topology on $M.$ This can be used to show that
if $G$ satisfies a ``strong F$\o$lner condition," then $H^1(G,M)$
is not Hausdorff.

As for the second result, we also show that when $\C G \subset M
\subset L^p(G)$ for some p $\in \N$, the $G$-module imbedding $\C
G \to M$ induces an imbedding of groups $H^1(G,\C G) \to
H^1(G,M).$ This (along with a minor result in [3]) shows that if
$H^1(G,M) = 0,$ then $G$ can have only one end.

The second paper (by E.G. Effros) characterizes the $\C$ algebra
$C^*_{red}(\Z^n)$ and provides a very explicit description of this
algebra. He shows first that $C^*_{red}(\Z^n) \cong C(\T^n)$ as
$\C$ algebras (where $\T^n$ is the $n$-torus). Showing that
$C^*_{red}(\Z^n)$ has no nontrivial idempotents of course then
shows that neither does $C(\T^n),$ which shows that $\T^n$ is
connected. Though this result is completely trivial, it ensues
some interesting mathematics and gives rise to the paper's title:
\textit{Why the Circle is Connected: An Introduction to Quantized
Topology}. What Effros's paper does mainly for ours is give us an
explicit isomorphism from $C(\T^n)$ onto $C^*_{red}(\Z^n)$.

Our first criterion for the Hausdorffness of $H^1(G,M)$ (Theorem
3) and the isomorphism between $C(\T^n)$ and $C^*_{red}(\Z^n)$
allow us to show that for $n \geq 1, H^1(\Z^n,C^*_{red}(\Z^n))$ is
not Hausdorff and so $\dim_{\C}H^1(\Z^n,C^*_{red}(\Z^n)) =
\infty.$

Secondly, the isomorphism from $C(\T^n)$ onto $C^*_{red}(\Z^n)$
allows us to describe the groups $H^n(\Z^n,C^*_{red}(\Z^n))$
explicitly. This in turn shows that there is a natural sequence of
imbeddings
$$0 \to H^1(\Z,C^*_{red}(\Z)) \to H^2(\Z^2,C^*_{red}(\Z^2)) \to
H^3(\Z^3,C^*_{red}(\Z^3)) \to \dots$$ This and the fact that
$H^1(\Z,C^*_{red}(\Z))$ is not Hausdorff allows us to show that
for each $n \geq 1,$ $\dim_{\C}H^n(\Z^n,C^*_{red}(\Z^n)) =
\infty.$

Finally, our second criterion for the Hausdorffness of $H^1(G,M)$
(Theorem 4) gives another (more simple) proof that
$H^1(\Z,C^*_{red}(\Z))$ is not Hausdorff.

I would like to give special thanks to my advisor, Dr. Peter
Linnell. His accessibility and the opinions that he gave me about
the problems I was working on were absolutely crucial to this
paper.

 \vspace{.1 in}

\noindent {\bf 2. Preliminaries and Definitions}

\vspace{.1 in}

\noindent Let $G$ be an infinite, finitely generated group.

\begin{definition} By a $G$ module $M$, we will mean a $\C$ vector
space $M$ along with a homomorphism of $G$ into $Aut(M).$
\end{definition}

\begin{definition} $G$ satisfies the strong F$\o$lner condition means
that for every finite subset $S$ of $G$ and every $\varepsilon
>0,$ there exists an $N \in \N$ such that for every $M \in \N$
there exists a finite subset $X$ of $G$ such that $|X|>M$ and
$|X-g\cdot X|<N$ for every $g \in S.$
\end{definition}

Because this definition may seem tedious, we should note that if
$G$ satisfies this condition, then $G$ must be amenable. Also, a
simple calculation shows that $\Z$ satisfies the strong F$\o$lner
condition.

\begin{definition} $\C G$ is the free $\C$ vector space with basis $G$.
In other words, $\C G$ is the set of all finite formal sums of the
form $\sum_g a_g g$ where the $a_g$ are in $\C$ and the $g$ are in
$G$.
\end{definition}

$\C G$ is a $G$ module in the obvious way ($G$ acts on an element
of $\C G$ by multiplication).

\begin{definition} For $1\leq p < \infty,$ $L^p(G)$ is the set of all formal sums (not
necessarily finite) of the form $\sum_g a_g g$ with the property
that $\sum_g |a_g|^p < \infty.$
\end{definition}

\begin{definition} The norm $||$ $||_p$ on $L^p(G)$ is given by
$$||\sum_g a_g g||_p = [\sum_g |a_g|^p]^{1/p}.$$
\end{definition}

$L^p(G)$ is complete in the topology induced by the norm $||$
$||_p$, so $L^p(G)$ is always a Banach space. There is a
multiplication defined on $\C G$ in the natural way ($a_i g_i
\cdot a_j g_j = (a_i a_j) (g_i g_j)$). While $\C G$ is a ring
(under componentwise addition and this multiplication as its ring
multiplication), $L^2(G)$ is not necessarily a ring. One can have
$\alpha,\beta \in L^2(G)$ with $\alpha\cdot\beta \notin L^2(G).$
However, it is well known that if $\alpha \in \C G$ and $\beta \in
L^2(G),$ then $\alpha\cdot\beta \in L^2(G).$ $L^2(G)$ is certainly
a Banach space over $\C$ and via this multiplication, $\C G$ may
be considered a subset of $B(L^2(G))$ (the bounded linear
operators on $L^2(G)).$

\begin{definition} The operator norm of an element $\alpha \in \C
G$ is given by $$||\alpha||_{op} = Sup
\{||\alpha\cdot\beta||_2:\beta \in L^2(G), ||\beta||_2 = 1\}.$$
\end{definition}

\begin{definition} $C^*_{red}(G)$ is the metric space completion
of $\C G$ under the operator norm $||$ $||_{op}$.
\end{definition}

We immediately have that $C^*_{red}(G)$ is composed entirely of
bounded linear operators on $L^2(G)$ and that $C^*_{red}(G)$ is
complete and therefore a Banach space. In addition, $\C G \subset
C^*_{red}(G) \subset L^2(G).$ We also know that $G$ acts on
$C^*_{red}(G)$ (in the obvious way) as continuous $\C$-linear
transformations.

\smallskip

Now we turn our attention to the groups $H^1(G,M).$ We will view
these groups in two different (though of course equivalent) ways.

First, let $A_G$ be the set of set maps $f\colon G \to M$ with the
property that for all $a$ and $b$ in $G,$ $a\cdot f(b) - f(ab) +
f(a) = 0.$ Let $B_G$ be the set of all such maps given by $f(b) =
b\cdot\alpha - \alpha$ for some fixed $\alpha \in M.$ Then we have
the following first definition of $H^1(G,M)$.

\begin{definition} $H^1(G,M)$ is the quotient group $A_G/B_G.$
\end{definition}

Now consider the ring $\Z G$ of all finite formal sums of the form
$\sum_g z_g g$ with $z_g \in \Z$ for all $g$ $\in G.$ Recall that
for an arbitrary ring $R$, an $R$ module is projective if and only
if it is a direct summand of a free $R$ module (there are several
definitions). We say that an infinite exact sequence of $\Z G$
modules $\dots \to P_1 \to P_0$ is a projective resolution of $\Z$
as $\Z G$ modules if it extends to an exact sequence $\dots \to
P_1 \to P_0 \to \Z \to 0$ of $\Z G$ modules. The $\Z G$ module
structure on $\Z$ is given by $g\cdot z_i$ = $z_i$ for $g \in G,
z_i \in \Z.$ Let $\dots \overset{d_2}{\to} P_1 \overset{d_1}{\to}
P_0 \overset{d_0}{\to} \Z \to 0$ be an exact sequence of $\Z G$
modules with $\dots \to P_1 \to P_0$ a projective resolution of
$\Z$. This first sequence induces another sequence $$0 \to Hom_{\Z
G}(\Z,M) \overset{d_0^*}{\to} Hom_{\Z G} (P_0,M)
\overset{d_1^*}{\to} Hom_{\Z G} (P_1,M)\overset{d_2^*}{\to}\dots$$
which gives us our second definition of $H^1(G,M)$.

\begin{definition} $H^1(G,M) = \ker (d^*_1)/\im (d^*_0).$
\end{definition}

As it turns out, this definition is independent of the choice of
projective resolution $\dots \longrightarrow P_1 \longrightarrow
P_0$ of $\Z$ as $\Z G$ modules. For the final part of our paper,
we will also need the following definition of $H^n(G,M)$.

\begin{definition} For n $\geq 1,$ $H^n(G,M) =
\ker(d^*_n)/\im(d^*_{n-1}).$
\end{definition}

Note that the two definitions of $H^1(G,M)$ are equivalent.

To define the topology on $H^1(G,M),$ we employ the first
definition of this group, $A_G/B_G.$ The topology is induced by
the topology of point-wise convergence on $A_G.$ That is for $f_n
\in A_G, f_n \longrightarrow 0$ means that $f_n(g) \in M$
converges to $0$ in the norm $||$ $||_M$ on $M$ for every $g \in
G.$ We should note that this is where the topologies on
$H^1(G,L^2(G))$ and on $H^1(G,M)$ differ. The basic open sets in
$A_G$ are the $f \in A_G$ such that $||f(g_1) -
a_1||_M,\dots,||f(g_n) - a_n||_M<\varepsilon$ for some choice of
fixed $g_i \in G, a_i \in M, \varepsilon>0$ and $n\in\N$. In other
words, the set of all such sets forms a basis for the topology on
$A_G$.

Finally, we will need to view $C(\T^n)$ as a subset of
$C(\T^{n+1})$ (especially in the proof of Theorem 8). We do this
as follows. For $f \in C(\T^n)$ and $(z_1,\dots,z_{n+1}) \in
\T^{n+1},$ define $f(z_1,\dots,z_{n+1}) = f(z_1,\dots,z_n).$

\vspace{.1 in}

\noindent {\bf 3. The group $H^1(G,M)$}

\vspace{.1 in}

\noindent Keep supposing that $G$ is infinite and finitely
generated. In addition, suppose that $M$ is a Banach space with
norm $||$ $||_M$ that is a left $\C G$ module and satisfies the
property that $G$ acts on $M$ as continuous $\C$-linear
transformations. Note that this implies that $M$ is a $G$ module.

\begin{theorem} Suppose that $\C G \subset M \subset L^p(G)$ for some
p $\in \N$. Then the $G$-module imbedding $\C G \to M$ induces an
imbedding of groups $H^1(G,\C G) \to H^1(G,M).$
\end{theorem}

\smallskip \noindent {\bf Proof:} Our result stated above has not
been published to date and could prove to be useful. We will,
however, follow the proof of Z.~Q.~Chen as described in [3]
(Proposition 1), which shows that the imbedding $\C G \to L^2(G)$
induces an imbedding $H^1(G,\C G) \to H^1(G,L^2(G)).$ Since $G$ is
finitely generated, suppose that $S$ is a finite generating set
for $G$. For an arbitrary $G$-module $N$, we define $C^n(G,N)$ to
be the set of all set maps from $G^n$ to $N$. In the case $n=0$,
we set $C^0(G,N) = N$. We have maps
$$d_0\colon C^0(G,\C G) \to C^1(G,\C
G),\hspace{.05 in}\delta_0\colon C^0(G,M) \to C^1(G,M)\hspace{.1
in} \textrm{and}$$
$$d_1\colon C^1(G,\C G) \to C^2(G,\C G),
\delta_1\colon C^1(G,M) \to C^2(G,M)$$ defined by $[d_0(a)](g) =
g\cdot a - a, a \in \C G, [\delta_0(f)](g) = g\cdot f - f, f \in
M, [d_1(h)](g_1,g_2) = g_1h(g_2) - h(g_1 g_2) + h(g_1),
[\delta_1(f)](g_1,g_2) = g_1f(g_2) - f(g_1 g_2) + f(g_1).$ Then
$$H^1(G,\C G) = \frac{\ker(d_1)}{\im(d_0)} \hspace{.1 in}
\textrm{and}\hspace{.1 in}H^1(G,M) =
\frac{\ker(\delta_1)}{\im(\delta_0)}.$$

The proposed imbedding is of course $a + (im(d_0)) \longmapsto a +
(im(\delta_0)),$ so define the natural homomorphism
$$\theta\colon \ker(d_1) \to
\frac{\ker(\delta_1)}{\im(\delta_0)}$$ via $a \longmapsto a +
(\im(\delta_0)).$ Since $\C G \subset M,$ we have that $\im(d_0)
\subset \im(\delta_0),$ so $\im(d_0) \subset \ker(\theta).$

Let $b\in \ker(\theta).$ Then $\im(b) \subset \C G$ and $b \in
\im(\delta_0).$ So there exists $f \in M$ such that for all $g \in
G, b(g) = g\cdot f - f \in \C G.$ We aim to show that $f$ must lie
in $\C G$ and thus $b \in im(d_0).$ Thus, we want to show that $f$
has finite support. Suppose $f = \sum a_g g.$ Note that for all $h
\in G, \hspace{.03 in} h\sum a_g g - \sum a_g g = \sum (a_g -
a_{hg})hg \in \C G.$ Thus, for for all $h \in G, \varphi(h) = \{g
\in G:a_g - a_{hg} \neq 0\}$ is finite. Since $S$ is a finite set,
it follows that
$$F(G) = \bigcup_{s \in S} \varphi(s)$$ is finite as well. We may
assume that for all $s \in S, s^{-1} \in S.$ Let $X$ be the Cayley
graph of $G$ with vertex set $G$ and edge set $\{(g,sg): s \in
S\}.$ By assuming that $S$ is closed under inverses, it follows
that $X$ can be viewed as an undirected graph. Thus, if $r$ and
$q$ (elements of $G$) are connected by an edge, it follows that
$r=sq$ for some $s \in S.$ Now consider the graph $X - F(G)$.
$F(G)$ is finite, so there are finitely many connected components
of $X-F(G)$. Then there exists a component of $X-F(G)$ that is
infinite (since $G$ is infinite). Let $q$ and $r$ be in this
connected component. So there exist $s_1,\dots,s_n \in S$ such
that $s_1\cdots s_n q = r$, and this path cannot pass through
$F(G)$, so we have the property that for all $i$, $s_i\cdots s_n q
\notin \varphi(s)$ for any $s \in S.$ Since $s_{i-1} \in S$ for
every $i$, it follows that $a_{s_{i-1}s_i\cdots s_n q} = a_{s_i
\cdots s_n q}$ for every $i$ and so $a_{q} = a_{s_n q} =
a_{s_{n-1} s_n q} = \dots = a_r.$ Thus, all of the $a_g$'s are
equal for all the $g$'s in $X-F(G)$. By virtue of there being
infinitely $g$'s in $X-F(G)$ and since $f = \sum_g a_g g \in M$
satisfies $\sum_g |a_g|^p < \infty,$ it follows that for all $g$
in this connected component of $X - F(G)$, $a_g = 0.$ Thus, if
$a_t \neq 0,$ it follows that $t$ lies in one of the finite
connected components of this graph, and there are only finitely
many such components (Since $F(G)$ is finite), so there are only
finitely many such $t$. i.e. $a_i = 0$ for all but finitely many
$i$, and thus $f \in \C G.$ $\Box$

\begin{lemma} $\dim_{\C}H^1(G,\C G) = |b(G)| - 1$ where $|b(G)|$ is
the number of ends of $G$.
\end{lemma}

\smallskip\noindent {\bf Proof:} For a proof of this, see [3](Lemma
2). $\Box$

\begin{theorem} Suppose that $\C G \subset M \subset L^p(G)$ for some p $\in \N.$ If $H^1(G,M) = 0,$ then $G$
has exactly 1 end.\end{theorem}

\smallskip\noindent {\bf Proof:} If $H^1(G,M) = 0,$ then by Theorem 1, $H^1(G,\C
G) = 0$. The previous lemma then tells us that $|b(G)| = 1$ ($G$
has exactly 1 end). $\Box$

\begin{theorem} $H^1(G,M)$ is not Hausdorff if and only if there
exists $f_i \in M$ with norm 1 ($||f_i||_M = 1$) for all $i$ with
the property that $||gf_i - f_i||_M \longrightarrow 0$ as $i
\longrightarrow \infty$ for every $g \in G.$
\end{theorem}

\smallskip \noindent {\bf Proof:} For this proof, we think of
$H^1(G,M)$ as the set $A_G$ of set maps $f$ from $G$ to $M$
satisfying $a\cdot f(b) - f(ab) + f(a) = 0$ for all $a,b$ $\in G$
modulo the set $B_G$ of maps $h$ of the form $h(g) = g\cdot e - e$
for some $e \in M.$ Saying that $H^1(G,M)$ is Hausdorff is
equivalent to saying that $0$ is closed in $H^1(G,M)$ which is of
course equivalent to saying that $B_G$ is closed in $A_G$. This is
equivalent to saying that $B_G$ is complete and is thus a Frechet
space. Remember that the topology on $A_G$ is that of point-wise
convergence, that is $f_n\colon G \to M$ tends to zero if and only
if $f_n(g) \longrightarrow 0 $ for every $g \in G$ in the norm
$||$ $||_M$ on $M$.

We have a continuous one-to-one map $h\colon M \to A_G$ whose
image is $B_G$ (the obvious map). Since we are no longer talking
about a topology on $B_G$ induced by $||$ $||_p,$ saying that this
map is continuous needs justification. Suppose that $(e_n) \in M$
satisfies $||e_n||_M \longrightarrow 0$. Fix $g \in G.$ Then
$||[h(e_n)](g)||_M = ||ge_n - e_n||_M \leq ||g-1||_M\cdot
||e_n||_M \longrightarrow 0$ as $n \longrightarrow \infty,$ since
$||g-1||_M$ is fixed. Since the topology on $B_G$ is that of
point-wise convergence, it follows that $h(e_n) \longrightarrow
0,$ so $h$ is continuous.

Now $M$ is a Banach space, so $M$ is certainly complete and a
Frechet space. Since a continuous, bijective map between Frechet
spaces has a continuous inverse, saying that $B_G$ is a Frechet
space is equivalent to saying that the inverse map from $B_G$ to
$M$ is continuous. We claim that this is equivalent to saying that
there does not exist a sequence $e_n \in M$ such that such that
$||e_n||_M = 1$ for all $n$ and $||ge_n - e_n||_M \longrightarrow
0$ for all $g \in G$.

The inverse map from $B_G$ to $M$ is given by the following. Say
$f\in B.$ Then $f$ is given by $f(g) = g\cdot e-e$ for some $e$
$\in M.$ This inverse map sends $f$ to $e$. Suppose that there
exists such a sequence $e_n \in M$ with $||e_n||_M = 1$ and
$||ge_n - e_n||_M \longrightarrow 0$ for all $g \in G.$ Thus, the
maps in $B_G$ determined by $e_n$ converge to $0$ (point-wise).
However, their image under this inverse map does not, which means
that this map cannot be continuous.

Conversely, suppose that this map is NOT continuous. Noting that
the topology on $B_G$ is that of point-wise convergence, this
means that there exists a sequence $a_n(g) = ge_n - e_n\in B_G$
such that $||a_n(g)||_M \longrightarrow 0$ for all $g \in G,$ but
$e_n$ does not converge to $0$ in $M.$ Then there exists
$\varepsilon > 0$ such that for all $N \geq 1,$ there exists $m
\geq N$ such that $||e_m||_M \geq \varepsilon.$ Multiply $e_n$ by
$1/\varepsilon$ and call the new sequence $e_n.$ Note that since
$\varepsilon > 0$ is fixed, the sequence of pre-images of the new
$e_n$ under this map still converge to $0$ for each $g \in G.$ Now
we can choose a subsequence $b_n$ of $e_n$ such that $||b_n||_M
\geq 1$ for all $n$. For each $m$, multiply $b_m$ by the unique
real number $0 < y_m \leq 1$ such that the new element of $M$
(rename it $X_m$) has norm $1$. Since $y_m \leq 1$ for all $m$, it
follows that for all $m$, $||X_m||_M = 1$ and $||gX_m -X_m||_M
\leq ||gb_m - b_m||_M \longrightarrow 0$ for all $g \in G.$ $\Box$

\begin{theorem} Keep the same assumptions as above about $G$ and
$M$ and suppose that $\C G \subset M \subset L^p(G)$ for some $p
\geq 2$. Suppose further that for every $\alpha \in \C G,
||\alpha||_1 \geq ||\alpha||_M \geq ||\alpha||_p.$ If $G$
satisfies the strong F$\o$lner condition, then $H^1(G,M)$ is not
Hausdorff.
\end{theorem}

\smallskip \noindent {\bf Proof:} Since $G$ is finitely generated,
we may let $G = \{g_1,g_2,\dots\}.$ For each $n\in\N,$ let $G_n =
\{g_1,\dots,g_n\}$. Because $G$ satisfies the strong F$\o$lner
condition, for each $n \in \N,$ we may pick $N_n \in \N$ such that
for every $M \in \N,$ there exists a finite subset $X$ of $G$ such
that $|X|>M$ and $|X- g_i \cdot X| < N_n$ for every $g_i \in S_n.$
Given any such $n$, choose a finite subset $X_n$ of $G$ such that
$|X_n|>(n\cdot N_n)^p$ and $|X_n - g_i \cdot X_n|<N_n$ for every
$g_i \in S_n.$ Let
$$\beta_n = \frac{\sum_{x \in X_n}x}{||\sum_{x \in X_n}x||_M}.$$
Note that for each $n$, $||\beta_n||_M = 1.$ Fix any $g_i \in G.$
Then for $n \geq i,$ we have
$$||g \cdot \beta_n - \beta_n||_M = \frac{||g\sum_{x \in X_n}
x - \sum_{x \in X_n}x||_M}{||\sum_{x \in X_n}x||_M}  \leq$$
$$\frac{||g\sum_{x \in X_n} x -
\sum_{x \in X_n}x||_M}{||\sum_{x \in X_n}x||_p} \leq
\frac{||g\sum_{x \in X_n} x - \sum_{x \in X_n}x||_1}{||\sum_{x \in
X_n}x||_p} =$$
$$=\frac{2\cdot |X_n - g_i\cdot X_n|}{||\sum_{x \in X_n}x||_p} <
\frac{2\cdot N_n}{||\sum_{x \in X_n}x||_p} = \frac{2\cdot
N_n}{(|X_n|)^{\frac{1}{p}}} < \frac{2\cdot N_n}{((n\cdot
N_n)^p)^{\frac{1}{p}}} = \frac{2}{n} \to 0.$$

Then $||\beta_n||_M = 1$ for all $n$ and $||g\cdot
\beta_n-\beta_n||_M \to 0$ for every $g \in G.$ By the previous
theorem, $H^1(G,M)$ is not Hausdorff. $\Box$

\begin{corollary} Suppose that $G$ satisfies the strong F$\o$lner
condition. Then $H^1(G,C^*_{red}(G))$ is not Hausdorff and
$\dim_{\C}(H^1(G,M)) = \infty.$ $\Box$
\end{corollary}

\vspace{.05 in}

\noindent {\bf 4. An Application: What Do the Groups
$H^1(\Z^n,C^*_{red}(\Z^n))$ and $H^n(\Z^n,C^*_{red}(\Z^n))$ look
like to $\C$?}

\vspace{.1 in}

Let $S(\T^n)$ be the square integrable functions on $\T^n.$ In
other words, $S(\T^n)$ is the set of functions $f\colon T^n\to \C$
such that the following integral exists and is finite:
$$\int_{|z_1| = 1}\cdots\int_{|z_n|=1} |f(z_1,\dots,z_n)|^2 dz_n\cdots dz_1.$$

Effros's paper ([2]) shows us that we have isomorphisms
$\theta_n\colon S(\T^n) \to L^2(\Z^n)$ and $\phi_n\colon
C(\T^n)\to C^*_{red}(\Z^n)$ both given by
$$f\longmapsto \sum_{m_1,\dots,m_n} a_{m_1,\dots,m_n} Z_1^{m_1}\cdots
Z_n^{m_n}\hspace{.1 in} \textrm{where}$$ $$a_{m_1,\dots,m_n} =
\int_{|z_1| = 1}\cdots\int_{|z_n|=1} Z_1^{-1-m_1} \cdots
Z_n^{-1-m_n}\cdot f(z_1,\dots,z_n)\cdot{1 \over (2\pi i)^n}\cdot
dz_n \cdots dz_1$$

\noindent For $f\in S(\T^n),$ Effros's paper ([2]) shows us that
we also have the following.
$$||\theta_n(f)||_2=[{1 \over (2\pi)^n}\int_{|z_1| = 1}\cdots\int_{|z_n|=1}
|f(z_1,\dots,z_n)|^2 dz_n \cdots dz_1]^{1/2}.$$

\begin{lemma}
For $f$ $\in C(\T^n),$ $||\phi_n(f)||_{op} = \sup_{z \in
\T^n}|f(z)|.$
\end{lemma}

\smallskip \noindent {\bf Proof:} Let $d$ be the usual distance metric on
$\C^n.$ First, suppose that $g$ $\in S(\T^n)$ satisfies
$||\theta_n(g)||_2 = 1.$
$$\Longrightarrow ||\phi_n(f)\theta_n(g)||_2=$$
$$[{1 \over (2\pi)^n}\int_{|z_1| = 1}
\cdots\int_{|z_n|=1}|g(z_1,\dots,z_n)|^2|f(z_1,\dots,z_n)|^2 dz_n
\cdots dz_1]^{1/2}=$$
$$\leq \sup_{z \in \T^n} |f(z)|\cdot [{1 \over (2\pi)^n}\int_{|z_1| = 1}
\cdots\int_{|z_n|=1}|g(z_1,\dots,z_n)|^2 dz_n \cdots
dz_1]^{1/2}=$$
$$=\sup_{z \in \T^n}|f(z)|.$$ $$\Longrightarrow
||\phi_n(f)||_{op}\leq \sup_{z \in \T^n} |f(z)|.$$

Now, without loss of generality, $||\phi_n(f)||_{op} = 1.$ Suppose
that $\sup_{z \in \T^n} |f(z)| > 1.$ So there is an $\varepsilon >
0$ such that $\sup_{z \in \T^n}|f(z)| > 1 + \varepsilon.$ Then
there is a $w$ $\in \T^n$ and $1 > \delta > 0$ such that $|f(z)| >
1 + \varepsilon$ for $d(z,w)$ $\leq \delta$ where $d$ is the usual
distance metric on $\C^n$ (This is due to the continuity of
$f(z)$). By $\int\cdots\int_{d(z,w)\leq \delta},$ we will mean
$\int\cdots\int_{\{d(z,w)\leq\delta\}\cap\T^n},$ the integral over
the intersection of the $\delta$ disk around $w$ (in $\C^n$) with
$\T^n$. Since
$$[{1 \over (2\pi)^n}\int_{|z_1| = 1}\cdots\int_{|z_n|=1}\cdot1\cdot
dz_n \cdots dz_1]^{1/2}= 1,$$ and since $\delta < 1,$ $\{z\in
\T^n:d(z,w)\leq \delta\}$ is properly contained in $\T^n.$ Thus,
there exists a $k \in \R^+$ such that $k > 1$ and a $g$ $\in
S(\T^n)$ defined by $g(z) = 0$ for $d(z,w)$ $> \delta$ and $g(z) =
k$ for $d(z,w)$ $\leq \delta$ with the property that
$$||\theta(g)||_2^2 = {1 \over (2\pi)^n}\int_{|z_1| = 1}\cdots
\int_{|z_n|=1}|g(z_1,\dots,z_n)|^2 dz_n\cdots dz_1 = 1.$$

Also, by the definition of $g$,
$$\int\cdots\int_{d(z,w)\leq \delta}|g(z_1,\dots,z_n)|^2 dz_n\cdots dz_1=$$
$$= \int_{|z_1| = 1}\cdots\int_{|z_n|=1} |g(z_1,\dots,z_n)|^2 dz_n \cdots dz_1.$$
Thus,
$$ 1 = ||\phi_n(f)||_{op} \geq ||\theta_n(g)\phi_n(f)||_2 =$$
$$= [{1 \over (2\pi)^n}\int_{|z_1| = 1}\cdots\int_{|z_n|=1}|g|^2|f|^2
dz_n\cdots dz_1]^{1/2} =$$
$$=[{1 \over (2\pi)^n}\int\cdots\int_{d(z,w)\leq \delta}|f|^2|g|^2
dz_n\cdots dz_1]^{1/2}=$$
$$\geq [{1 \over (2\pi)^n}\int\cdots\int_{d(z,w)\leq \delta}(1 +
\varepsilon)^2|g|^2 dz_n\cdots dz_1]^{1/2} =$$
$$= (1 + \varepsilon)[{1 \over (2\pi)^n}\int\cdots\int_{d(z,w)\leq
\delta}|g|^2 dz_n\cdots dz_1]^{1/2}=$$
$$= (1 + \varepsilon)[{1 \over (2\pi)^n}\int_{|z_1|=
1}\cdots\int_{|z_n|=1}|g|^2 dz_n\cdots dz_1]^{1/2} = 1 +
\varepsilon > 1,$$ A contradiction. Thus, $||\phi_n(f)||_{op} =
\sup_{z \in \T^n}|f(z)|.$ $\Box$

\begin{theorem} Let $n \geq 1.$ Then $H^1(\Z^n,C^*_{red}(\Z^n))$
is not Hausdorff and so

\noindent $\dim_{\C}H^1(\Z^n,C^*_{red}(\Z^n)) = \infty.$
\end{theorem}

\smallskip \noindent {\bf Proof:} Again, let $d$ be the usual
distance metric on $\C^n$. Suppose that
$H^1(\Z^n,C^*_{red}(\Z^n))$ is Hausdorff. By Theorem $3$, there
cannot exist $(f_i) \in C^*_{red}(\Z^n) = C(\T^n)$ such that
$||\phi_n(f_i)||_{op} = 1$ for all $i$ with the property that
$||\phi_n(gf_i - f_i)||_{op} \longrightarrow 0$ as $i
\longrightarrow \infty$ for every $g$ $\in \Z^n.$ We will exhibit
such an $(f_i).$

Since given an $i$, $\{(z_1,\dots,z_n) \in
\T^n:d((z_1,\dots,z_n),(1,\dots,1))>1/i \}$ and
$\{(z_1,\dots,z_n):d((z_1,\dots,z_n),(1,\dots,1)) < 1/(2i)$ for
all j$\}$ are disjoint open subsets of $\T^n$ with positive
distance between them (this distance is $1/(2i)$), for each $i$ we
have a function $f_i \in C(\T^n)$ with image contained in the
closed unit disk that is identically $0$ on the first set and
identically $1$ on the second. By the previous lemma,
$||\phi_n(f_i)||_{op} = \sup_{z \in \T^n}|f_i(z)| = 1$ for all
$i$. Choose $(m_1,\dots,m_n) \in \Z^n$. Let $B_i$ =
$\{(z_1,\dots,z_n) \in \T^n:d((z_1,\dots,z_n),(1,\dots,1))\leq 1/i
\}$ Then since each $f_i$ is zero on the set $A_i =
\{(z_1,\dots,z_n) \in \T^n:d((z_1,\dots,z_n),(1,\dots,1))>1/i \}$,
$$||\phi_n(Z_1^{m_1}\cdots Z_n^{m_n}f_i(z_1,\dots,z_n) - f_i(z_1,\dots,z_n))||_{op} =$$
$$= \sup_{(z_1,\dots,z_n) \in \T^n} |Z_1^{m_1}\cdots
Z_n^{m_n}f_i(z_1,\dots,z_n) - f_i(z_1,\dots,z_n)|=$$
$$= \sup_{B_i}|Z_1^{m_1}\cdots
Z_n^{m_n}f_i(z_1,\dots,z_n) - f_i(z_1,\dots,z_n)|\leq$$
$$\leq \sup_{B_i}|Z_1^{m_1}\cdots
Z_n^{m_n} - 1|\cdot Sup_{B_i} |f_i|=$$
$$= \sup_{B_i}|Z_1^{m_1}\cdots
Z_n^{m_n} - 1| \longrightarrow 0 \hspace{.1 in}\textrm{as }
\hspace{.1 in} i \longrightarrow \infty.$$

This is due to the fact that on $B_i,
d((Z_1,\dots,Z_n),(1,\dots,1))\leq 1/i$ and the fact that the
$m_1,...,m_n$ are fixed. In other words, $||gf_i - f_i||_{op}
\longrightarrow 0$ as $i \longrightarrow \infty$ for every $g \in
\Z^n.$ $\Box$

Now we turn our attention to the groups
$H^n(\Z^n,C^*_{red}(\Z^n)).$ We view $\T^n$ as
$\{(z_1,\dots,z_n):|z_i| = 1$ for all $ i\}.$ Again from Effros's
paper ([2]), we know that for all $n$, $C^*_{red}(\Z^n) \cong
C(\T^n)$ as $\Z^n$ modules. But what is the $\Z^n$ module
structure on $C(\T^n)?$ We view $\Z^n$ as the free abelian group
on the generators $x_1,\dots,x_n.$ Then the action of an element
$x_1^{e_1}\cdots x_n^{e_n}$ on an element $f(z_1,\dots,z_n)$ of
$C(\T^n)$ is given trivially by $(x_1^{e_1}\cdots x_n^{e_n})\cdot
f(z_1,\dots,z_n) = (z_1^{e_1}\cdots z_n^{e_n})f(z_1,\dots,z_n).$
First, we must recall a case of the Kunneth theorem for projective
resolutions.

\begin{theorem} (\emph{K\"unneth Theorem}) Let $G$ be a group and let $P\colon\dots \to P_1
\overset{d_1}{\to} P_0, Q\colon \dots\to Q_1
\overset{\delta_1}{\to} Q_0$ be projective resolutions of $\Z$ as
$\Z G$ modules. Then $P\otimes Q$ is a projective resolution of
$\Z$ as $\Z G$ modules where $[P\otimes Q]_n = \oplus_{i+j=n}
P_i\otimes Q_j.$ The boundary maps $\partial_n$ are given by
$\partial_n(p_i\otimes q_j) = d_i(p_i)\otimes q_j + (-1)^j
p_i\otimes \delta_j(q_j)$ on simple tensors $p_i\otimes q_j.$
\end{theorem}

\smallskip \noindent {\bf Proof:} This is Theorem V.2.1 of [4].
$\Box$

We aim to give an explicit description of the groups
$H^n(\Z^n,C^*_{red}(\Z^n))$ and to show that there is a natural
sequence of imbeddings
$$0 \to H^1(\Z,C^*_{red}(\Z)) \to H^2(\Z^2,C^*_{red}(\Z^2)) \to
H^3(\Z^3,C^*_{red}(\Z^3)) \to \dots$$ In order to do so, we must
prove the following lemma.

\begin{lemma} Let $G=\Z^n$ be the free abelian group on $X_1,\dots, X_n.$
Then there exists a projective resolution of $\Z$ as $\Z G$
modules with exactly n+1 nonzero modules (including $\Z$): $0 \to
\Z G \overset{\rho_n}{\to} (\Z G)^n \to \dots \to \Z \to 0$ where
the map $\rho_n$ is given by $\rho_n(1) = [(X_1 - 1),-(X_2 -
1),\dots,-(X_n - 1)].$

\end{lemma}

\smallskip \noindent {\bf Proof:} Note that for all $n$, $\Z G =
\Z[X_1,\dots,X_n,X_1^{-1},\dots,X_n^{-1}]$. For $n = 1$, we know
that there is a projective resolution of $\Z\colon 0 \to \Z G
\overset{d}{\to} \Z G \to \Z \to 0,$ where $d$ is given by
$d(f(X_1)) = (X_1 - 1)f(X_1).$ Suppose that $n$ is at least $2$.

Next we must note that for any $n$, $\Z G \cong \Z \Z \otimes_{\Z}
\Z \Z \otimes_{\Z} \dots \otimes_{\Z} \Z \Z$ where there are $n$
terms, and the isomorphism is an isomorphism of $\Z G$ modules.
The isomorphism $\theta\colon \Z G \to \Z \Z \otimes_{\Z} \dots
\otimes_{\Z} \Z \Z$ is given explicitly by $1$ $\longmapsto
1\otimes \dots \otimes 1.$ The action of $G$ (the free abelian
group on $X_1,\dots,X_n$) on $\Z \Z \otimes_{\mathbb{Z}} \dots
\otimes_{\mathbb{Z}} \Z \Z$ is as follows:
$$X_1^{r_1}\cdots X_n^{r_n} * a_1 \otimes \dots \otimes a_n =
X_1^{r_1}a_1 \otimes \dots \otimes X_n^{r_n}a_n.$$

Let us look again at the projective resolution of $\Z$ as a $\Z G$
module for $n = 1$ given by $0 \to \Z G \overset{d}{\to} \Z G
\overset{aug}{\to} \Z \to 0.$ In order to properly view this
resolution as a resolution, we will call the first $\Z
G$\hspace{.03 in} $P_1$ and the second $P_0.$ Thus, we have a
resolution $0 \to P_1 \overset{d}{\to} P_0 \overset{aug}{\to} \Z
\to 0,$ where $P_1 = P_0 = \Z \Z.$

The K\"unneth theorem tells us that the tensor product of
resolutions of $\Z$ is again a resolution, so by induction on $n$,
for all $n$ at least $2$, there is a resolution of $\Z$:
$$0 \to (P_1 \otimes_{\Z} \dots \otimes_{\Z} P_1)
\overset{d_n}{\to} (P_1 \otimes_{\Z} P_1 \otimes_{\Z} \dots
\otimes_{\Z} P_0)\oplus \dots \oplus (P_0\otimes_{\Z} P_1
\otimes_{\Z} \dots \otimes_{\Z} P_1) \to$$
$$\dots \to \Z \to 0$$

where everywhere $\dots$ appears implies that there are exactly
$n$ objects being either $\oplus$ed or $\otimes$ed together,
except the last $\dots$ which implies that there are exactly $n +
1$ (including $\Z$) nonzero terms in the sequence. The map $d_n$
is given explicitly from the definition of tensor products of
resolutions by
$$d_n(a_1 \otimes \dots \otimes a_n) = $$
$$[(X_1 - 1)a_n \otimes a_2 \otimes \dots \otimes a_n, -(a_1 \otimes
(X_2 - 1)a_2 \otimes \dots \otimes a_n),\dots,-(a_1 \otimes \dots
\otimes a_{n-1} \otimes (X_n -1)a_n)].$$

Note that the map $d_n$ is just a generalization of our original
map $d:$ $P_1 \to P_0.$ We know that
$$(P_1 \otimes_{\Z} \dots \otimes_{\Z} P_1 \otimes_{\Z} P_0) \oplus
\dots \oplus(P_0 \otimes_{\Z} P_1 \otimes_{\Z} \dots \otimes_{\Z}
P_1) \cong (\Z G)^n$$ via a map $\zeta$ where the isomorphism is
as $\Z G$ modules (and where $G$ is still the free abelian group
on $X_1,\dots,X_n$). The map $\zeta$ sends $d_n(1 \otimes \dots
\otimes 1) = [((X_1 - 1)\otimes 1 \otimes \dots \otimes 1),
-(1\otimes (X_2 - 1) \otimes \dots \otimes 1),- \dots, -(1 \otimes
\dots \otimes 1 \otimes (X_n -1))]$ \hspace{.03 in} to \hspace{.03
in} $[X_1 - 1, -(X_2 - 1), \dots ,-(X_n - 1)]$. Translating back
into $\Z G$ lingo, it follows directly from the K\"unneth Theorem
that we have a projective resolution of $\Z$ as $\Z G$ modules
with exactly $n+1$ (including $\Z$) nonzero modules
$$0 \to \Z G \overset{\rho_n}{\to} (\Z
G)^n \to \dots \to \Z \to 0.$$

The map $\rho_n: \Z G \to (\Z G)^n$ is then of course given by a
composition of maps:
$$\rho_n(1) = \zeta(d_n(\theta(1))) = \zeta(d_n(1\otimes \dots
\otimes 1)) = [X_1 - 1, -(X_2 - 1), \dots , - (X_n - 1)].
\hspace{.1 in} \Box$$

Now we are able to give an explicit description of the groups
$H^n(\Z^n,C^*_{red}(\Z^n)).$

\begin{theorem} Let $G = \Z^n$ be the free abelian group on $X_1, \dots ,X_n.$
Then we have that

$$H^n(G,C^*_{red}(G)) \cong \frac{C(\T^n)}{(Z_1 - 1)C(\T^n) + \dots + (Z_n -
1)C(\T^n)}.$$

\end{theorem}

\smallskip\noindent {\bf Proof:} Since for all $n$, $C^*_{red}(\Z^n) \cong
C(\T^n),$ it suffices to prove the statement for $H^n(G,C(\T^n))$
in the place of $H^n(G,C^*_{red}(G)).$ From the previous lemma, we
have a projective resolution of $\Z$ with $n+1$ nonzero terms
(including $\Z$):
$$0 \to \Z G \overset{\rho_n}{\to} (\Z
G)^n \to \dots \to \Z \to 0.$$ From this, we get a new sequence
$B$:
$$0 \to Hom_{\Z G}(\Z,C(\T^n)) \to \dots \to Hom_{\Z G}((\Z G)^n,C(\T^n))
\overset{\rho_n^*}{\to} Hom_{\Z G}(\Z G,C(\T^n))\to 0.$$

Using this sequence, we may compute the cohomology group
$$H^n(G,C^*_{red}(G)) = H^n(G,C(\T^n)) = \frac{Hom_{\Z G}(\Z
G,C(\T^n))}{\rho_n^*(Hom_{\Z G}((\Z G)^n,C(\T^n)))}.$$

Define $\theta\colon Hom_{\Z G}(\Z G,C(\T^n))\to C(\T^n)/[(Z_1 -
1)C(\T^n) + \dots + (Z_n - 1)C(\T^n)]$ by $\phi \longmapsto
\phi(1) + [(Z_1 - 1)C(\T^n) + \dots + (Z_n - 1)C(\T^n)].$ This map
is certainly onto, because if $F \in C(\T^n),$ we may just define
$\sigma \in Hom_{\Z G}(\Z G,C(\T^n))$ by $\sigma(1) = F.$ Now we
can compute $ker(\theta).$

Suppose that $\phi \in Hom_{\Z G}((\Z G)^n,C(\T^n)).$ Then
$$\phi(\rho_n(1)) = \phi(X_1 - 1, - (X_2 - 1),\dots, -(X_n - 1)) =$$
$$(X_1 - 1)*\phi(1,0,\dots,0) + \dots + (X_n - 1)*\phi(0,\dots,0,-1)\in$$
$$ \in [(Z_1 - 1)C(\T^n) + \dots + (Z_n - 1)C(\T^n)].$$
$$\Longrightarrow \rho_n^*(Hom_{\Z G}((\Z G)^n,C(\T^n))) \subset
ker(\theta).$$

Now suppose that $\phi \in ker(\theta).$ So $\phi(1) = (Z_1 -
1)f_1 + \dots + (Z_n - 1)f_n$ for some $f_i \in C(\T^n).$ Define
$\pi \in Hom_{\Z G}((\Z G)^n,C(\T^n))$ by $\pi(1,0,\dots,0) = f_1,
\pi(0,1,0,\dots,0) = -f_2, \dots, \pi(0,\dots,0,1) = -f_n.$ Then
$\pi(\rho_n(1)) = \phi(1)$, and therefore we have that $\phi =
\pi(\rho_n) \in \rho_n^*(Hom_{\Z G}((\Z G)^n,C(\T^n))).$
$$\Longrightarrow H^n(G,C^*_{red}(G)) = H^n(G,C(\T^n)) = $$
$$\frac{Hom_{\Z G}(\Z G,C(\T^n))}{\rho_n^*(Hom_{\Z G}((\Z
G)^n,C(\T^n)))} \cong \frac{C(\T^n)}{(Z_1 - 1)C(\T^n) + \dots +
(Z_n - 1)C(\T^n)}. \hspace{.1 in} \Box$$

\begin{theorem} There is a sequence of group (and $\C$ algebra)
imbeddings

$$0 \to H^1(\Z,C^*_{red}(\Z))
\overset{d_1}{\to} H^2(\Z^2,C^*_{red}(\Z^2)) \overset{d_2}{\to}
H^3(\Z^3,C^*_{red}(\Z^3)) \to \dots$$

\end{theorem}

\smallskip \noindent {\bf Proof:} By the previous theorem, for all
$n \geq 1,$
$$H^n(\Z^n,C^*_{red}(\Z^n)) = \frac{C(\T^n)}{(z_1 - 1)C(\T^n) + \dots + (z_n
- 1)C(\T^n)}.$$ Fix an $n$ and define
$$\phi\colon \frac{C(\T^n)}{(z_1 - 1)C(\T^n) + \dots + (z_n -
1)C(\T^n)}$$
$$\to \frac{C(\T^{n+1})}{(z_1 - 1)C(\T^{n+1}) + \dots +
(z_{n+1} - 1)C(\T^{n+1})}$$
$$\textrm{via}\hspace{.3 in} f + [(z_1 - 1)C(\T^n) + \dots + (z_n -
1)C(\T^n)]$$ $$ \longmapsto f + [(z_1 - 1)C(\T^{n+1}) + \dots +
(z_{n+1} - 1)C(\T^{n+1})].$$

Since $(z_1 - 1)C(\T^n) + \dots + (z_n - 1)C(\T^n) \subset (z_1 -
1)C(\T^{n+1}) + \dots + (z_{n+1} - 1)C(\T^{n+1}),$ $\phi$ is
certainly well defined. To see that $\phi$ is injective, we just
let $n=1$ (the other cases are identical). If $f$ $\in C(\T)$ and
$f(z_1) = (z_1 - 1)f_1(z_1,z_2) + (z_2 - 1)f_2(z_1,z_2)$ with $
f_1,f_2 \in C(\T^2),$ then since $f$ doesn't vary at all with
$z_2,$ it follows that we can set $z_2 =1$ and then we have that
$f(z_1) = (z_1 - 1)f(z_1,1) \in (z_1 - 1)C(\T).$ Thus, $\ker(\phi)
= 0$ and so $\phi$ is injective. $\Box$

\begin{theorem} For $n \geq 1,$

$$\dim_{\C}H^n(\Z^n,C^*_{red}(\Z^n)) = \infty.$$

\end{theorem}

\smallskip \noindent {\bf Proof:} For all $n$, by the last theorem
$H^1(\Z,C^*_{red}(\Z))$ imbeds in $H^n(\Z^n,C^*_{red}(\Z^n)).$
Since $\dim_{\C}H^1(\Z,C^*_{red}(\Z)) = \infty,$ the result
follows. $\Box$

\smallskip

Finally, using Theorem 4, we are able give a more simple proof of
the following result (one proof is given above in Theorem 5).

\begin{theorem} $H^1(\Z,C^*_{red}(\Z))$ is not Hausdorff and
$\dim_{\C}H^1(\Z,C^*_{red}(\Z)) = \infty.$
\end{theorem}

\smallskip \noindent {\bf Proof:} We know that for $\alpha \in \C
\Z, ||\alpha||_1 \geq ||\alpha|_{op} \geq ||\alpha||_2.$ Thus,
Theorem 4 tells us that $H^1(\Z,C^*_{red}(\Z))$ is not Hausdorff.
$\Box$

\vspace{.1 in}

\noindent {\bf 5. Our Results (Summarized)}

\vspace{.1 in}

\noindent Let $G$ be an infinite, finitely generated group. Let
$M$ be a Banach space (with norm $||$ $||_M$) that is also a left
$\C G$ module such that $G$ acts on $M$ as continuous $\C$-linear
transformations. In summary, the following useful results were
proved in this paper and do not appear in any papers to date.

\vspace{.05 in}

\textbf{(1)} Suppose that $\C G \subset M \subset L^p(G)$ for some
p $\in \N$. Then the $G$-module imbedding $\C G \to M$ induces an
imbedding of groups $H^1(G,\C G) \to H^1(G,M).$

\textbf{(2)} If $\C G \subset M \subset L^p(G)$ for some p $\in
\N$ and if $H^1(G,M) = 0,$ then $G$ has exactly $1$ end.

\textbf{(3)} $H^1(G,M)$ is not Hausdorff if and only if there
exists $f_i \in M$ with norm 1 ($||f_i||_M = 1$) for all $i$ with
the property that $||gf_i - f_i||_M \longrightarrow 0$ as $i
\longrightarrow \infty$ for every $g \in G.$

\textbf{(4)} If $\C G \subset M \subset L^p(G)$ for some $p \geq
2$, $||\alpha||_1\geq ||\alpha||_M \geq ||\alpha||_p$ for every
$\alpha \in \C G$, and if $G$ satisfies the strong F$\o$lner
condition, then $H^1(G,M)$ is not Hausdorff and therefore
$\dim_{\C}H^1(G,M) = \infty.$

\vspace{.05 in}

We have used these results (with our motivation coming from
Effros's paper ([2])) to show the following for $n \geq 1$.

\vspace{.05 in}

\textbf{(1)} $H^1(\Z^n,C^*_{red}(\Z^n))$ is not Hausdorff and thus
$\dim_{\C}H^1(\Z^n,C^*_{red}(\Z^n)) = \infty$ (with two proofs for
$n = 1$).

\textbf{(2)} $\dim_{\C}H^n(\Z^n,C^*_{red}(\Z^n)) = \infty.$

\textbf{(3)} There is a natural sequence of imbeddings
$$0 \to H^1(\Z,C^*_{red}(\Z)) \to H^2(\Z^2,C^*_{red}(\Z^2)) \to
H^3(\Z^3,C^*_{red}(\Z^3)) \to \dots$$

\textbf{(4)} $$H^1(\Z^n,C^*_{red}(\Z^n)) = \frac{C(\T^n)}{(Z_1 -
1)C(\T^n) + \dots + (Z_n - 1)C(\T^n)}.$$

\vspace{.15 in}

\noindent {\bf 6. Conjectures and Future Work}

\begin{conjecture} Let $G$ be an infinite, finitely generated group. Let $M$ be a
Banach space (with norm $||$ $||_M$) that is also a left $\C G$
module such that $G$ acts on $M$ as continuous $\C$-linear
transformations. Then $G$ is amenable if and only if $H^1(G,M)$ is
not Hausdorff.
\end{conjecture}

We should note that the proof given in [3] of this with $L^2(G)$
is specific to that case.

\end{document}